# Flexibility Management for Space Logistics via Decision Rules

Hao Chen[1]
*Georgia Institute of Technology, Atlanta, GA, 30332*

Brian M. Gardner[2] and Paul T. Grogan[3]
*Stevens Institute of Technology, Hoboken, NJ, 07030*

and
Koki Ho[4]
*Georgia Institute of Technology, Atlanta, GA, 30332*

**This paper develops a flexibility management framework for space logistics mission planning under uncertainty through decision rules and multi-stage stochastic programming. It aims to add built-in flexibility to space architectures in the phase of early-stage mission planning. The proposed framework integrates the decision rule formulation into a network-based space logistics optimization formulation model. It can output a series of decision rules and generate a Pareto front between the expected mission cost (i.e., initial mass in low-Earth orbit) and the expected mission performance (i.e., effective crew operating time) considering the uncertainty in the environment and mission demands. The generated decision rules and the Pareto front plot can help decision-makers create implementable policies immediately when uncertainty events occur during space missions. An example mission case study about space station resupply under rocket launch delay uncertainty is established to demonstrate the value of the proposed framework.**

## Nomenclature

$\mathcal{A}$ = arc set

$b$ = remaining stock, kg

$c$ = cost matrix

---

[1] Ph.D. Student, Department of Aerospace Engineering, AIAA Student Member.
[2] Ph.D. Student, School of Systems and Enterprises, AIAA Student Member.
[3] Assistant Professor, School of Systems and Enterprises, AIAA Senior Member.
[4] Assistant Professor, Department of Aerospace Engineering, AIAA Member.

| | | |
|---|---|---|
| $D$ | = | delay length, days |
| $\mathbf{d}$ | = | demand or supply vector, kg |
| $\mathcal{G}$ | = | time-expanded network |
| $H$ | = | concurrency constraint matrix |
| $h$ | = | crew time loss, days |
| $J$ | = | mission cost function |
| $\mathcal{K}$ | = | scenario set |
| $\mathcal{L}$ | = | rocket launch set |
| $M$ | = | large constant number |
| $\mathcal{N}$ | = | node set |
| $p$ | = | scenario probability |
| $Q$ | = | commodity transformation matrix |
| $R$ | = | safety stock decision rule |
| $r$ | = | available safety stock, kg |
| $\mathcal{T}$ | = | time step set |
| $u$ | = | additional supply, kg |
| $\mathcal{V}$ | = | spacecraft set |
| $\mathbf{W}$ | = | time window vector |
| $x$ | = | commodity flow variable, kg or # |
| $\mathcal{Z}$ | = | mission performance |
| $\alpha$ | = | safety stock availability (binary) |
| $\beta$ | = | additional supply necessity (binary) |
| $\gamma$ | = | weighting coefficient |
| $\varepsilon$ | = | total number of commodities, # |
| $\theta$ | = | total number of concurrency constraints, # |
| $\eta$ | = | commodity consumption rate, kg/day |
| $\psi$ | = | commodity shortage penalty, day/kg |
| $\zeta$ | = | rocket launch time function |

*Subscripts*

| | | |
|---|---|---|
| *i* | = | node index |
| *j* | = | node index |
| *k* | = | scenario index |
| *l* | = | launch index |
| *t* | = | time step index |
| *v* | = | spacecraft index |

## I. Introduction

THE complexity of space exploration programs increases with the increasing interests of human and robotic space exploration. In traditional space programs like Apollo, each mission was logistically independent. Thus, the uncertainties in one mission (e.g., spacecraft flight delays/failures, mission demand changes) had a relatively limited impact on subsequent missions, which simplified the evaluation of consequences and solutions to handle uncertainties. More recent studies show that sustainable space exploration needs to be accomplished by a campaign composed of multiple interdependent missions [1-3]. Leveraging a campaign enables effective space exploration with a combination of different logistics paradigms such as pre-deployment, carry-along, and resupply. At the same time, however, mission interdependencies also introduce the potential of "cascading failure" behaviors in space exploration, similar to those in the interconnected infrastructures in terrestrial civil infrastructure systems [4]. Particularly, a rocket launch delay in an earlier mission may significantly impact the performance of later missions, such as further launch delays in the future and supply shortages to support mission operations. Therefore, there is a growing need to consider the uncertainties in launch delay in the campaign design.

One example of the impact of complex interdependencies between space missions can be found in the International Space Station (ISS) program. The ISS has encountered three unexpected mission failures in an 8-month span from October 2014: two caused by rocket launch failures [i.e., Orbital Commercial Resupply Service (CRS)-3 and SpaceX CRS-7] and one caused by spacecraft flight failure (i.e., Progress 59) [5]. These failures resulted in an immediate loss in excess of 6,832 kg cargo and supplies including the first of two docking adapters to support the commercial crew program [6-8]. They also led to delays in subsequent space missions to identify the reason for the failures and incorporate the necessary changes to support the return to flight certification for the launch vehicle or spacecraft. The return to flight time was 720 days for the Orbital Antares launch vehicle and 285 days for the SpaceX Falcon launch vehicle [5]. Since the ISS program utilized multiple cargo suppliers and maintained a level of supply to ensure astronaut safety in the event of an extended disruption, the impact on ISS operations was minimal and there was no risk to astronaut safety [9]. However, future large-scale space campaigns are more complex than ISS. They can include multiple space stations and habitat resupply missions in parallel. To mitigate the impact of the uncertain environment on mission operations, we need a systematic methodology to quantitatively identify the level of safety stock balancing the mission cost and mission performance while ensuring the safety of astronauts under uncertain mission environments.

In the space logistics research field, multiple frameworks have been proposed in the past to perform efficient space mission planning under deterministic environments. These frameworks were established through heuristic methods [10], graph theory [11], simulations [12], and multi-commodity network flow models [1-3]. However, deterministic methods tend to give overly optimistic designs and bias anticipated mission performance. Space mission operations following the mission planning results optimized under deterministic environments may end up with failure or significant mission cost increase if any uncertain event occurs during the space campaign such as spacecraft flight delays.

The research for design under uncertainty in the space field mainly focused on analyzing the impact of uncertainty on space missions, but little work has addressed how to counter these uncertainties in the mission planning phase. For example, Shull [13] discussed the impact of campaign level risks, such as flight delay and cancellation, in a human lunar exploration architecture. Added flexibility was integrated with the format of pre-positioning stockpile and evaluated through sensitivity analysis. Moreover, the threat of uncertainty for future human Mars missions was studied by Stromgren et al. [14] using the Exploration Maintainability Analysis Tool (EMAT). However, these studies cannot generate directly implementable decision strategies for decision-makers to follow in response to uncertain events.

Beyond the space research area, various studies have considered the flexibility in design from different perspectives, including taxonomy [15], empirical evaluation [16], and stochastic programming with real options analysis [17] or decision rules [18]. However, these proposed methods cannot be applied directly to our problem of space logistics planning under uncertainties in launch delays. First, mathematically, the transportation of additional safety stock needs to be taken into account during the logistics of other pre-determined mission demands because of the limited time window opportunities (e.g., due to the launch environment and orbital dynamics). The decision rule formulation needs to be integrated into the space logistics optimization problem. Furthermore, there is a conceptual difference between our problem and classical real-option problems in that our problem involves the uncertainties of the time windows, which prevents the decision makers from responding to the uncertainties in real time. In classical built-in flexibility examples, we typically assume that we can respond to the uncertainties when they realize (e.g., expanding the building in response to the increasing demand). However, in our problem, the time steps of the decision windows themselves are uncertain. Therefore, when there is a launch delay, the impact appears on the mission performance immediately, but decision-makers cannot react to it until the delayed launch itself (i.e., after the impacts have already occurred). Namely, at the point of decision making, we can only prepare for the future and compensate for the uncertain event immediately before it, rather than reacting to the uncertain event in real time.

In response to this background, this paper proposes a decision rule-based flexibility management framework for space logistics mission design leveraging multi-stage stochastic programming. It aims to add built-in flexibility to space architectures in the phase of early-stage mission planning. Decision rules, also called implementable policies, map the observations of uncertainty data to the decisions directly. Thus, the resulting decision rule can be implemented as general decision strategies for decision-makers to follow when stochastic events occur.

In this research, space logistics optimization problems are formulated based on network flow models, and decision rules are formulated focusing on the uncertainties of rocket launch delays. However, the proposed decision rule formulation can be easily extended to integrate with other space logistics optimization methods and handling different types of uncertainty sources. Leveraging our previous work [19], a case study on crewed space station resupply missions is established to demonstrate how flexible space mission designs achieve greater mission performance under uncertainty as compared to mission concepts developed under deterministic assumptions.

Our work provides an important step to make space mission planning and architecture design flexible to counter potential uncertainties in space missions. The proposed methodology is particularly useful for evaluating the performance of space infrastructure design under stochastic mission environments and developing a reliable space mission schedule for future campaign-level space exploration.

The remainder of this paper is organized as follows. Section II first introduces the problem setting and proposes the formulation of decision rules for space logistics mission planning problems. The performance of the proposed flexibility management framework is then demonstrated in Sec. III through a space station resupply mission example. Finally, Sec. IV concludes the paper and discusses future works.

## II. Methodology

**A. Problem Setting**

In this paper, we consider a problem of space logistics mission planning under uncertain rocket launch delay. The logistics requirement is to support persistent scientific missions/experiments on space stations by regular rocket launches. However, the rocket launch may delay, which would impact the operation of space missions. We assume that the space missions/experiments to be supplied are time-sensitive so that the supply shortage influences the operating time of these missions within the planned mission period. The total length of the mission period is fixed regardless of the launch delay, which means the mission duration cannot be extended even if it is temporarily paused due to the shortage of scientific instruments and maintenance spares. Thus, a temporary supply shortage reduces the time that astronauts can spend on the experiments or maintenance, resulting in less available operating time to support the mission performance. Note that this assumption is not always true for any space missions. For space missions with a specific time window, the mission can be non-extendable; whereas some space experiments may be extended with no or little penalty, such as educational activities. In reality, the penalty may also be in the format of extra cost instead of the operating time loss. In this paper, we presume this assumption to simplify the analysis process of the mission planning and the structure of decision rules.

To mitigate such operating time loss, we can launch more supplies in earlier launches and maintain them as safety stocks for later missions. These safety stocks are buffers to temporarily support space activities. Launching more safety stocks can decrease the operating time loss but increase mission costs in earlier space missions. Therefore, the decision maker's goal is a tradeoff between the mission cost and performance (i.e., operating time loss)[5]. The amount of safety stock for each rocket launch is the outcome of decisions. Note that the focus of this paper is on the scientific experiment cost and performance; the crew consumables or other mandatory supplies are assumed to be sufficiently stocked on the stations and thus are not considered as part of the tradeoff.

---

[5] Crew time also has been used as a metric of space missions in some previous literature such as SpaceNet [12].

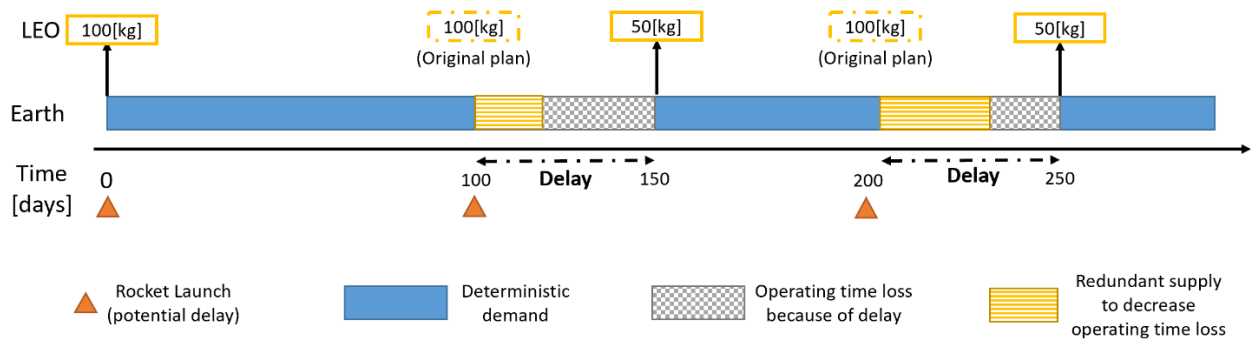

**Fig. 1 How delay and safety stock influence space missions (an example).**

Fig. 1 shows one example of how the rocket launch delay and safety stock influence space missions in our problem. Consider a case where there is a rocket launch every 100 days to ISS. For this example, we assume ISS needs materials to support an in-space experiment which is consumed at a rate of 1 kg/day. If all supply missions are launched as planned, each rocket will deliver 100 kg of material to ISS, which is enough to support a 100-day duration space mission. If there is a 50-day delay in the second and the third launches as shown in Fig. 1, we only need to launch 50 kg of material to support the rest space experiment in each launch. However, because of the shortage of material on ISS, the space experiment is paused during the 50-day delays. If there is no safety stock on ISS, two 50-day delays result in the suspension of the space experiment for 100 days due to the lack of supply. Our space mission is still terminated on the originally scheduled end date. As a result, the actual time for the experiment, which is originally planned as 300 days long, only runs for 200 days because of the shortage of materials. This results in 100 days of operating time loss due to the crew not having the material available to run the science experiment for the full 300 days duration. If we launched more material in the first and the second rocket launches, the safety stock can partially support the space experiments and decrease the operating time loss when a rocket launch delay occurs in subsequent missions. However, the mission costs for the first and the second rocket launch increase to launch these additional safety stock materials.

## B. Scenario Generation Model

The scenario generation model provides launch dates for the cargo and crew to support mission requirements as uncertainty inputs. The launch delay model is prescribed as a continuous cumulative probability function created from available data. To limit the number of possible mission outcomes to a computationally manageable number, we perform inverse transform sampling to generate uncertainty mission scenarios based on the continuous launch delay probability function. In the later case study, we consider the case where each launch delay is independent, but the approach can be applied to dependent launch delays as well. Because we conduct a random sampling method, the probability of realization for each possible outcome is equivalent. If there are $n$ scenarios in total, then the probability is $1/n$ for each mission scenario.

**C. Flexibility Management Background and Overview**

To manage space mission operations under complex uncertain mission scenarios, we need to incorporate uncertainty mitigation methods in space logistics optimization. Multiple analytical tools have been introduced in past literature to deal with uncertainties in engineering fields, such as the backward induction and dynamic programming by Buurman et al. [20], the optimization-based hybrid real options analysis by Jiao [21], and the simulation-based decision analysis by de Neufville et al. [22]. However, these methods provide decisions as outcomes of the model. There is no general decision strategy for managers to follow when a different scenario appears throughout the system life cycle. On the other hand, the decision rule, also called implementable policy, maps the observations of uncertainty data to the decisions. Cardin et al. [18] showed that for a highly structured problem, optimization based on decision rules can provide a good approximation of optimal solution that is typically found by real option analysis techniques, while the solutions of decision rule were more practical and intuitive for decision-makers to use in practical applications. Numerical experiments on waste-to-energy systems in past literature showed that decision rules-based methods can find a very close solution to traditional uncertainty analysis methods with a significantly reduced computational time required [18]. These results also showed that decision rules are suitable to analyze complex systems with multiple sources of uncertainties. Multiple decision strategies can also be considered concurrently.

The decision rule is a gradient-free direct policy search method. Compared with the value-based method in reinforcement learning, direct policy search solves the problem by searching the policy space directly without the value estimation of states or state-action pairs. For certain problems, including the space logistics problem considered in this paper, where we have a clear idea about what the nature of policy representations might look like, constructing good policies is often easier than modeling good state functions. Another benefit of the decision rule method is that it does not rely on gradient information. It can be explicitly integrated with the space logistics mission planning framework through mixed-integer linear programming (MILP), where a global optimum can be achieved. Compared with other stochastic sequential decision-making approaches, considering uncertainties via decision rules also avoids the dramatic increase of computational complexity.

There are multiple types of decision rules that are suitable for different uncertainty problems. Four classes of decision rules were identified and studied by Garstka and Wets via case examples [23], including zero-order, linear, safety-first, and conditional-go decision rules. A detailed introduction about decision rules is provided in Ref. [18] and [23]. Conditional-go decision rule has been used in generation expansion planning problems [24]; whereas linear decision rule has been used in electrical reserve operation problems [25]. In this paper, we focus on the uncertainty from rocket launch delays, which is contrary to traditional uncertainty in other engineering fields [18, 23-25]. When a rocket launch delay occurs, we have already missed the opportunity to respond to this situation. We have to make decisions in previous launches to counter the uncertainty before we observe the delay.

To identify the influence of launch delay to mission objectives and determine the amount of safety stock for each launch, the decision rule in rocket launch delay problems is formulated based on a combination of the conditional-go decision rule and the linear decision rule through integer programming (IP). The amount of safety stock is the decision variable. The conditional-go decision rule identifies whether previous safety stock on board is enough to counter the launch delay and the penalty of launch delay to the mission performance. The linear decision rule finds the additional supply needed to be launched in the upcoming mission. The developed decision rule can also be used in other uncertainties in space logistics problems, such as spacecraft flight or docking delay and demand change in space missions.

To balance mission cost and mission performance for flexibility management, we consider a multi-objective optimization solved through the weighted sum method. The initial mass in low-Earth orbit (IMLEO) and crew operating time loss are used as the metrics for mission cost and mission performance. IMLEO is widely used as a space logistics mission cost metric in previous studies, such as space logistics mission planning frameworks [1-3] and the Mars mission design reference architecture by NASA [26]. The crew operating time loss, as discussed earlier, is used to represent space mission performance. More specifically, the operating time loss is calculated as the summation of the crew time loss for each commodity, representing its dependence on both the length of launch delays and the crew time task assignment. For example, astronauts may spend a different length of time on scientific experiments and maintenance. The consumption rate of scientific instruments and maintenance spares may also be different. Therefore, for the same length of supply shortage, the lack of scientific instruments or maintenance spares can cause different length of operating time loss.

The flexibility management framework can be expressed as the following optimization problem.

Minimize

$$\sum_{k \in \mathcal{K}} p_k \mathcal{J}_k(\pmb{x}_k) + \gamma \sum_{k \in \mathcal{K}} p_k \mathcal{Z}_k(\pmb{R}) \qquad (1)$$

Subject to

$$f(\pmb{x}_k, \pmb{R}, \pmb{y}) \leq 0 \qquad (2)$$

$$g(\pmb{x}_k, \pmb{R}, \pmb{y}) = 0 \qquad (3)$$

where $\mathcal{J}_k$ is the mission cost function for mission scenario $k$ from the space logistics problem; $\mathcal{Z}_k$ represents the operating time loss for mission scenario $k$ due to the lack of supplies determined from the decision rule formulation. The decision variables $\pmb{x}_k$ and $\pmb{R}$ denotes the space transportation commodity flow and the safety stock for space logistics and decision rule, respectively. Note that $\pmb{x}_k$ and $\pmb{R}$ are dependent on each other through the constraints, and so they need to be optimized concurrently. In the constraints, $\pmb{y}$ represents other variables in mission planning that do not directly impact the value of objectives. In this problem, the mission cost $\mathcal{J}_k$, space logistics decision variable $\pmb{x}_k$, and the mission performance $\mathcal{Z}_k$ all depend on the mission scenario $k$, while the decision rule variable $\pmb{R}$ is independent of mission scenarios as it is a general implementable policy that can be applied to any scenario. We define a mission scenario set $\mathcal{K}$, then the probability for each scenario is $p_k = 1/\|\mathcal{K}\|$.

In Eq. (1), $\gamma$ is the weighting coefficient to determine how important the mission performance is compared with the mission cost. If $\gamma$ is infinity, decision-makers will prepare for the worse scenario, where safety stock levels reach maximum for each launch. By changing the weighting coefficient $\gamma$, we can get Pareto front solutions between the expected IMLEO $\mathcal{J}$ and the expected operating time loss $\mathcal{Z}$ based on different strategies. Note that we claim the solutions as the Pareto front since all the decision rule solutions obtained by the weighted sum method are all on the Pareto-optimal front. However, the weighted sum method cannot guarantee to find all Pareto-optimal solutions, particularly in the case of a nonconvex objective space. We use the weighted sum method because the objective function keeps linear that guarantees a global optimum of the optimization. Moreover, it also avoids extra constraints to ensure computational efficiency. The constraints Eqs. (2) and (3) are general constraint expressions for space logistics and decision rules that will be introduced in detail in the following sections.

In the following sections, we first introduce the space logistics optimization model to calculate $\mathcal{J}$ based on mission demands and decision rules. Then, we introduce the proposed decision rule method to determine $\mathcal{Z}$ based on the stochastic mission operation environment. Finally, we integrate decision rules into space logistics mission planning through the optimization problem Eq. (1) and establish the flexibility management framework.

### D. Mission Cost Evaluation through Space Logistics

This section introduces the space logistics model to calculate the mission cost (i.e., IMLEO) $J$. The space logistics section in the flexibility management framework is developed based on network-based multi-commodity flow models [1-3]. It was originally proposed to optimize spacecraft design, space infrastructure design, and space transportation scheduling concurrently under a deterministic environment. This network-based method models orbits or planets as nodes; transportation trajectories as arcs. Crew, propellant, spacecraft, and all payloads are considered as commodities flowing along arcs. Consider a time-expanded network $\mathcal{G}$, which is made up of a set of arcs $\mathcal{A}$. The arc set $\mathcal{A}$ is defined based on a set of nodes $\mathcal{N}$ (index $i, j$), a set of vehicles $\mathcal{V}$ (index $v$), and a set of time steps $\mathcal{T}$ (index $t$). We also need to define a set of scenarios $\mathcal{K}$ (index $k$), where each scenario $k$ has a probability of $p_k$ (i.e., $\sum_{k \in \mathcal{K}} p_k = 1$) to occur. The decision variable for space logistics optimization is the commodity flow variable $\boldsymbol{x}_{kvijt}$, representing the commodity flow from node $i$ to node $j$ using spacecraft $v$ at time $t$ in the scenario $k$. For each scenario $k$, we need a demand parameter $\boldsymbol{d}_{kit}$, representing demands or supplies of different commodities in node $i$ at time $t$. Mission demands are negative and mission supplies are positive. To measure the mission cost, a commodity cost coefficient $\boldsymbol{c}_{vijt}$ is defined. If there are $\varepsilon$ type of commodities in total, $\boldsymbol{x}_{kvijt}$, $\boldsymbol{d}_{kit}$, and $\boldsymbol{c}_{vijt}$ are all $\varepsilon \times 1$ vectors.

Besides the pre-planned mission demands and supplies defined in $\boldsymbol{d}_{kit}$, we also need another demand variable $\boldsymbol{u}_{kit}$, representing additional commodities required to be delivered to maintain the level of required safety stocks for each rocket launch. It is an interface variable to connect space logistics optimization and decision rules. The value of this variable is positive only at the node of space stations. We define a subset of the node set $\mathcal{N}$ for space stations, $\widehat{\mathcal{N}}$. Therefore, we have

$$\boldsymbol{u}_{kit} = \begin{cases} \boldsymbol{0}_{\varepsilon \times 1} \text{ if } i \in \mathcal{N} \setminus \widehat{\mathcal{N}} \\ \boldsymbol{u}_{kit}, \text{otherwise} \end{cases}$$

where the value of $\boldsymbol{u}_{kit}$, when $i \in \widehat{\mathcal{N}}$, is determined by the decision rule for each space station.

Besides the parameters defined above, we also need the following parameters for the space logistics optimization formulation:

$\Delta t_{ij}$ = time of flight.

$Q_{vij}$ = commodity transformation matrix.

$H_{vij}$ = concurrency constraint matrix.

$W_{kij}$ = mission time windows.

Then, the formulation for space logistics optimization under uncertainties is shown as follows.

Minimize

$$\sum_{k \in \mathcal{K}} p_k \mathcal{J}_k = \sum_{k \in \mathcal{K}} p_k \left( \sum_{(v,i,j,t) \in \mathcal{A}} \left( \boldsymbol{c}_{vijt}^T \boldsymbol{x}_{kvijt} \right) \right) \quad (4)$$

Subject to

$$\sum_{(v,j):(v,i,j,t) \in \mathcal{A}} \boldsymbol{x}_{kvijt} - \sum_{(v,j):(v,j,i,t) \in \mathcal{A}} Q_{vji} \boldsymbol{x}_{kvji(t-\Delta t_{ji})} \leq \boldsymbol{d}_{kit} - \boldsymbol{u}_{kit} \quad \forall i \in \mathcal{N} \ \forall t \in \mathcal{T} \ \forall k \in \mathcal{K} \quad (5)$$

$$H_{vij} \boldsymbol{x}_{kvijt} \leq \boldsymbol{0}_{\theta \times 1} \quad \forall (v,i,j,t) \in \mathcal{A} \ \forall k \in \mathcal{K} \quad (6)$$

$$\begin{cases} \boldsymbol{x}_{kvijt} \geq \boldsymbol{0}_{\varepsilon \times 1} & \text{if } t \in W_{kij} \\ \boldsymbol{x}_{kvijt} = \boldsymbol{0}_{\varepsilon \times 1} & \text{otherwise} \end{cases} \quad \forall (v,i,j,t) \in \mathcal{A} \ \forall k \in \mathcal{K} \quad (7)$$

$$\boldsymbol{x}_{kvijt} \in \mathbb{R}_{\geq 0}^{\varepsilon \times 1} \quad \forall (v,i,j,t) \in \mathcal{A} \ \forall k \in \mathcal{K}$$

In this formulation, Eq. (4) is the objective function, which calculates the expected mission cost, where $p_k = 1/\|\mathcal{K}\|$. When evaluating the expected IMLEO, the cost coefficients for the launch arcs are set as one, while the rest are set as zero. Equation (5) is the mass balance constraint that limits commodity flows to satisfy the demands. The second term $Q_{vji} \boldsymbol{x}_{kvji(t-\Delta t_{ji})}$ represents commodity transformations during space transportation and space mission operations, including propellant burning, crew consumables, and resource generations by space infrastructures. The additional supply $\boldsymbol{u}_{kit}$ is the interface variable connecting space logistics and decision rules. Equation (6) is the spacecraft concurrency constraint, which defines the upper bound of commodity flows limited by the spacecraft propellant capacity and the payload capacity. The number of concurrency constraints considered in mission planning is denoted by $\theta$. Equation (7) is the time window constraint defined by the time window vector $W_{kij}$. Only when time windows are open, are space flights permitted.

### E. Operating Time Loss Evaluation through Decision Rules

This section introduces the decision rule method to determine the mission performance metric (i.e., crew operating time loss) $\mathcal{Z}$. For a flexibility management problem for multiple space stations, decision rules need to be established independently for each space station. For simplicity, in this section, we only show the decision rule formulation for one space station. The decision rules for multiple space stations and the combination with space logistics will be discussed in the next section. Thus, we omit the node index in the additional supply variable and rewrite it as $\boldsymbol{u}_{kt}$.

For a space station, we define a set of rocket launches $\mathcal{L}$ (index $l$). The safety stock, denoted by $\boldsymbol{R}_l$ for the rocket launch $l$, is the decision rule variable to be optimized. The decision rule can be expressed as "launch up to $\boldsymbol{R}_l$ kg safety stock in rocket launch $l-1$". After observing the delay in rocket launch $l-1$, when the system is available to launch a rocket again, decision-makers prepare safety stock for the next mission. Now, since launch delay $l-1$ has already occurred, decision-makers know the amount of remaining stock $\boldsymbol{b}_{k(l-1)}$ at the destination, where $k$ is the index of uncertainty scenario. Based on the decision rule variable $\boldsymbol{R}_l$, which is the safety stock we have to maintain for launch $l$, we can determine the additional supply $\boldsymbol{u}_{kt}$ that should be launched in addition to the necessary mission demand in launch $l-1$, which is $\boldsymbol{R}_l - \boldsymbol{b}_{k(l-1)}$. However, there is a possibility that the remaining stock $\boldsymbol{b}_{k(l-1)}$ is larger than $\boldsymbol{R}_l$, where no additional supply is needed. Therefore, the expression of additional supply is $\boldsymbol{u}_{kt} = \max(\boldsymbol{R}_l - \boldsymbol{b}_{k(l-1)}, \boldsymbol{0}_{\varepsilon \times 1})$, where $t$ is the time index in space mission planning formulation. This expression works when $t$ is the time step at which launch $l-1$ happens. Then, for each mission scenario $k$, we know the available safety stock prepared for launch $l$, defined as $\boldsymbol{r}_{kl} = \max(\boldsymbol{R}_l, \boldsymbol{b}_{k(l-1)})$. The available safety stock $\boldsymbol{r}_{kl}$ can also be written as $\boldsymbol{r}_{kl} = \boldsymbol{u}_{kt} + \boldsymbol{b}_{k(l-1)}$. Based on the available safety stock for each scenario $k$, we can calculate the length of the supply shortage $\boldsymbol{h}_{kl}$ for each launch delay, in the unit of days. In this research, we measure mission performance through the total operating time loss $\mathcal{Z}$, which can be then calculated based on $\boldsymbol{h}_{kl}$. The optimization framework eventually outputs $\boldsymbol{R}_l$ as the decision rule result. Decision-makers can make decisions directly and immediately after a launch delay occurs based on the value of $\boldsymbol{R}_l$.

Besides the aforementioned notations, we also need to define the following parameters for the decision rule formulation.

$\boldsymbol{c}''$ = operating time loss weighting coefficient for different commodities.

$D_{kl}$ = length of the launch delay, days.

$\boldsymbol{\eta}$ = commodity consumption rate, kg/day.

$\boldsymbol{\psi}$ = commodity shortage penalty.

Then, the decision rule optimization can be formulated as follows.

Minimize
$$\sum_{k \in \mathcal{K}} p_k \mathcal{Z}_k = \sum_{k \in \mathcal{K}} p_k \sum_{l \in \mathcal{L}} \boldsymbol{c}''^T \boldsymbol{h}_{kl} \tag{8}$$

Subject to
$$\begin{cases} \boldsymbol{u}_{kt} = \max(\boldsymbol{R}_l - \boldsymbol{b}_{k(l-1)}, \boldsymbol{0}_{\varepsilon \times 1}) & \text{if } t = \zeta(l-1) \\ \boldsymbol{u}_{kt} = \boldsymbol{0}_{\varepsilon \times 1} & \text{otherwise} \end{cases} \quad \forall t \in \mathcal{T} \quad \forall l \in \mathcal{L} \quad \forall k \in \mathcal{K} \tag{9}$$

$$\boldsymbol{h}_{kl} = \max\left((D_{kl} \boldsymbol{1}_{\varepsilon \times 1} - (\boldsymbol{u}_{kt} + \boldsymbol{b}_{k(l-1)})) \oslash \boldsymbol{\eta}\right) \circ \boldsymbol{\psi}, \boldsymbol{0}_{\varepsilon \times 1}) \quad \forall l \in \mathcal{L} \quad \forall k \in \mathcal{K} \tag{10}$$

$$\boldsymbol{b}_{kl} = \max((\boldsymbol{u}_{kt} + \boldsymbol{b}_{k(l-1)}) - D_{kl}\boldsymbol{\eta}, \boldsymbol{0}_{\varepsilon \times 1}) \quad \forall l \in \mathcal{L} \ \forall k \in \mathcal{K} \tag{11}$$

$$\boldsymbol{R}_l, \boldsymbol{u}_{kt}, \boldsymbol{h}_{kl}, \boldsymbol{b}_{kl} \in \mathbb{R}_{\geq 0}^{\varepsilon \times 1} \quad \forall l \in \mathcal{L} \ \forall k \in \mathcal{K}$$

Equation (8) is the objective function, which calculates the expected operating time loss, where $p_k = 1/\|\mathcal{K}\|$. Equation (9) determines the value of the interface variable $\boldsymbol{u}_{kt}$. It connects the additional supply $\boldsymbol{u}_{kt}$ in launch $l$ and the remaining stock $\boldsymbol{b}_{k(l-1)}$ after launch $l-1$ with the expected safety stock $\boldsymbol{R}_l$ for mission $l$. The rocket launch time index function, $t = \zeta(l-1)$, outputs the time step of each rocket launch. This function converts the time index $t$ in the space logistics formulation into the rocket launch index $l$ in the decision rule section. Equation (10) calculates the operating time loss because of the shortage of certain commodities. Note that, $\oslash$ and $\circ$ are Hadamard operators. They are the division and product of each element with the same indices between two matrices. These two matrices should have the same dimension. If the safety stock is sufficient to support the space mission during delays, $D_{kl}\boldsymbol{1}_{\varepsilon \times 1} - (\boldsymbol{u}_{kt} + \boldsymbol{b}_{k(l-1)}) \oslash \boldsymbol{\eta}$ is negative, which then enforces the operating time loss to be zero. Equation (11) calculates the level of remaining stock. If the safety stock is not sufficient to support a launch delay, the remaining supplies $\boldsymbol{b}_{kl}$ are zeros. But if there are enough available safety stock $\boldsymbol{r}_{kl} = \boldsymbol{u}_{kt} + \boldsymbol{b}_{k(l-1)}$ to support the mission during a launch delay, there are some remaining supplies after this delay. The remaining commodities become part of the safety stock for the next launch.

Note that, even though the constraints in the above decision rule optimization formulation are all equality constraints, they are all nonlinear functions. To solve this problem effectively, we can linearize the constraints by introducing two binary variable vectors $\boldsymbol{\alpha}_{kl}$ and $\boldsymbol{\beta}_{kl}$ to represent whether there is enough safety stock to support the mission operation and whether additional supplies are needed for a specific rocket launch mission, respectively.

$$\boldsymbol{\alpha}_{kl} = \begin{cases} 0 \text{ if } D_{kl} \leq (\boldsymbol{u}_{kt} + \boldsymbol{b}_{k(l-1)}) \oslash \boldsymbol{\eta} \\ 1 \text{ if } D_{kl} > (\boldsymbol{u}_{kt} + \boldsymbol{b}_{k(l-1)}) \oslash \boldsymbol{\eta} \end{cases}$$

$$\boldsymbol{\beta}_{kl} = \begin{cases} 0 \text{ if } \boldsymbol{b}_{k(l-1)} \leq \boldsymbol{R}_l \\ 1 \text{ if } \boldsymbol{b}_{k(l-1)} > \boldsymbol{R}_l \end{cases}$$

Define a large constant scalar $M$, the decision rule formulation can be linearized as follows.

Minimize

$$\sum_{k \in \mathcal{K}} p_k Z_k = \sum_{k \in \mathcal{K}} p_k \sum_{l \in \mathcal{L}} \boldsymbol{c}''^T \boldsymbol{h}_{kl} \tag{12}$$

Subject to

$$\begin{cases} \begin{cases} \boldsymbol{u}_{kt} \geq \boldsymbol{R}_l - \boldsymbol{b}_{k(l-1)} \\ \boldsymbol{u}_{kt} \leq M\boldsymbol{\beta}_{kl} + \boldsymbol{R}_l - \boldsymbol{b}_{k(l-1)} \\ \boldsymbol{u}_{kt} \leq M(\boldsymbol{1}_{\varepsilon\times 1} - \boldsymbol{\beta}_{kl}) \end{cases} \text{if } t = \zeta(l-1) \\ \boldsymbol{u}_{kt} = \boldsymbol{0}_{\varepsilon\times 1} \quad \text{otherwise} \end{cases} \quad \forall t \in \mathcal{T} \quad \forall l \in \mathcal{L} \quad \forall k \in \mathcal{K} \quad (13)$$

$$\begin{cases} \boldsymbol{h}_{kl} \geq \left(D_{kl}\boldsymbol{1}_{\varepsilon\times 1} - (\boldsymbol{u}_{kt} + \boldsymbol{b}_{k(l-1)}) \oslash \boldsymbol{\eta}\right) \circ \boldsymbol{\psi} \\ \boldsymbol{h}_{kl} \leq M(\boldsymbol{1}_{\varepsilon\times 1} - \boldsymbol{\alpha}_{kl}) + \left(D_{kl}\boldsymbol{1}_{\varepsilon\times 1} - (\boldsymbol{u}_{kt} + \boldsymbol{b}_{k(l-1)}) \oslash \boldsymbol{\eta}\right) \circ \boldsymbol{\psi} \quad \forall l \in \mathcal{L} \quad \forall t = \zeta(l-1) \quad \forall k \in \mathcal{K} \quad (14) \\ \boldsymbol{h}_{kl} \leq M\boldsymbol{\alpha}_{kl} \end{cases}$$

$$\begin{cases} \boldsymbol{b}_{kl} \geq (\boldsymbol{u}_{kt} + \boldsymbol{b}_{k(l-1)}) - D_{kl}\boldsymbol{\eta} \\ \boldsymbol{b}_{kl} \leq M\boldsymbol{\alpha}_{kl} + (\boldsymbol{u}_{kt} + \boldsymbol{b}_{k(l-1)}) - D_{kl}\boldsymbol{\eta} \quad \forall l \in \mathcal{L} \quad \forall t = \zeta(l-1) \quad \forall k \in \mathcal{K} \quad (15) \\ \boldsymbol{b}_{kl} \leq M(\boldsymbol{1}_{\varepsilon\times 1} - \boldsymbol{\alpha}_{kl}) \end{cases}$$

$$\begin{cases} D_{kl}\boldsymbol{1}_{\varepsilon\times 1} - (\boldsymbol{u}_{kt} + \boldsymbol{b}_{k(l-1)}) \oslash \boldsymbol{\eta} \geq M(\boldsymbol{\alpha}_{kl} - \boldsymbol{1}_{\varepsilon\times 1}) \\ D_{kl}\boldsymbol{1}_{\varepsilon\times 1} - (\boldsymbol{u}_{kt} + \boldsymbol{b}_{k(l-1)}) \oslash \boldsymbol{\eta} \leq M\boldsymbol{\alpha}_{kl} \end{cases} \quad \forall l \in \mathcal{L} \quad \forall t = \zeta(l-1) \quad \forall k \in \mathcal{K} \quad (16)$$

$$\begin{cases} \boldsymbol{b}_{k(l-1)} - \boldsymbol{R}_l \geq M(\boldsymbol{\beta}_{kl} - \boldsymbol{1}_{\varepsilon\times 1}) \\ \boldsymbol{b}_{k(l-1)} - \boldsymbol{R}_l \leq M\boldsymbol{\beta}_{kl} \end{cases} \quad \forall l \in \mathcal{L} \quad \forall k \in \mathcal{K} \quad (17)$$

$$\boldsymbol{R}_l, \boldsymbol{u}_{kt}, \boldsymbol{h}_{kl}, \boldsymbol{b}_{kl} \in \mathbb{R}_{\geq 0}^{\varepsilon\times 1} \quad \boldsymbol{\alpha}_{kl}, \boldsymbol{\beta}_{kl} \in \{0,1\}^{\varepsilon\times 1} \quad \forall l \in \mathcal{L} \quad \forall k \in \mathcal{K}$$

In this formulation [Eqs. (12)-(17)], all constraints become linear constraints. They are equivalently linearized constraints of the original nonlinear equality constraints. The additional constraints Eqs. (16) and (17) determines the value of binary variables $\boldsymbol{\alpha}_{kl}$ and $\boldsymbol{\beta}_{kl}$.

### F. Flexibility Management Framework for Space Logistics

Substituting the formulations to calculate $\mathcal{J}$ and $\mathcal{Z}$ proposed in Sec. III.D and Sec. III.E into Eq. (1), we combine decision rule with space logistics and obtain a flexibility management framework for space logistics mission design under uncertainties, shown as follows.

Minimize

$$\sum_{k\in\mathcal{K}} p_k \mathcal{J}_k + \gamma \sum_{k\in\mathcal{K}} p_k \mathcal{Z}_k = \sum_{k\in\mathcal{K}} p_k \left(\sum_{(v,i,j,t)\in\mathcal{A}} \left(\boldsymbol{c}_{vijt}^T \boldsymbol{x}_{kvijt}\right)\right) + \gamma \sum_{k\in\mathcal{K}} p_k \sum_{l\in\mathcal{L}} \boldsymbol{c}''^T \boldsymbol{h}_{kl} \quad (18)$$

Subject to

*Constraint set 1 (space logistics model)*

$$\text{Eqs. (5)-(7)}$$

*Constraint set 2 (decision rules)*

*For each space station $i \in \widehat{\mathcal{N}}$:*

$$\text{Eqs. (13)-(17)}$$

Equation (18) is the objective function. The first term calculates the value of expected mission cost $\mathcal{J}$ (i.e., IMLEO) and the second term calculates the expected value of mission performance $\mathcal{Z}$ (i.e., operating time loss). Note that the decision rules represented by Eqs. (13-17) in constraint set 2 are defined for each space station $i \in \widehat{\mathcal{N}}$. Thus, the corresponding decision variables also need to be defined separately for each station.

This framework models the space logistics problem considering uncertainties as a MILP formulation. It is a single objective optimization formulation with weighted sum objectives. The flexibility management framework introduced in this section is suitable to be directly used for mission design under uncertainties in rocket launch delay and spacecraft flight delay. It can also be extended to consider other traditional uncertainty factors in space logistics problems.

In the next section, the performance of this flexibility management framework for space logistics is evaluated based on a space station resupply mission case study.

## III. Case Study: Space Station Resupply Logistics

**A. Space Station Resupply Mission Background**

In this paper, the performance of the flexibility management framework is evaluated based on a crewed space station resupply example mission in a lunar orbit. The mission scenario is developed based on the data of past resupply missions for ISS and the planned Gateway concept [27].

*1. Mission Demand*

The example mission planning case study is a one-year multi-mission logistics to the space station. The case study considers a crew of four with four missions per year: two missions with cargo and crew and two missions with cargo only. A crew rotation will be performed during each crew mission. Each group of crew stays on the space station for six months. The goal of flexibility management is to satisfy the cargo demand for space activities while guaranteeing astronaut safety during space station operations. The estimated yearly cargo needs for the space station is up to 16,750 kg based on ISS pressurized cargo needs defined in the CRS2 Request for information [28]. The crew consumable requirement is defined as 17.1 kg/day for 4 crew based on NASA Design Reference Architecture 5.0 [26]. The remaining cargo is split between science and maintenance cargo. A review of the recent NASA CRS mission overviews (6 December 2015 through 13 June 2017) shows a nearly even split between science (37.2%) and maintenance cargo mass (36.4%) of the total cargo [7, 29-37], detailed in Table 1. For this example mission case study, it is expected that the initial maintenance needs of the considered lunar orbital station will be lower than the ISS assuming that it is significantly smaller than the ISS. Therefore, the prescribed cargo split is 66% for science and 33% for maintenance for the example scenario. The consumption rate for the science and maintenance cargo is assumed to be evenly distributed over the mission duration of one year.

We also need to set penalties for the shortage of science instruments and maintenance spares to measure the influences of rocket launch delays and determine the operating time loss. This penalty is relevant to the importance of commodities to the mission operation. We assume astronauts will spend 80% of their time performing scientific experiments and 20% of their time performing maintenance on the space station. Thus, one day shortage of scientific instruments will lead to a 0.8 day operating time loss, while one day shortage of maintenance spares will lead to a 0.2 day operating time loss. We will perform a sensitivity analysis to evaluate the impact of the crew operating time assignment later. Based on these assumptions, our goal is to find the optimal amount of safety stock to reduce operating time loss while balancing the mission cost at the same time.

**Table 1 ISS CRS Cargo Manifest Mass by Category (Dec. 2015 to Jun. 2017)**

| Category | Mass, kg | Percentage of total, % |
|---|---|---|
| Crew Supplies | 4,893 | 26.4 |
| Science Investigations | 6,881 | 37.2 |
| Other | 6,746 | 36.4 |
| Total (Including Packaging) | 19,334.1 | 100 |

The orbit of the space station is considered to be the Near Rectilinear Halo Orbit (NRHO). It takes about 5 days to fly from low Earth orbit (LEO) or back to LEO. The total $\Delta V$ from LEO to NRHO is 3.53 km/s, while the $\Delta V$ from NRHO to LEO is 3.51 km/s [38].

Based on the consumption rates of science instruments, maintenance spares, and crew consumables introduced above, the demands and supplies of this space station resupply mission are listed in Table 2. The mission plan is shown in Fig. 2.

This paper focuses on the decision rule optimization and analysis. Therefore, for simplicity, in this case study, no spacecraft design is considered and an existing spacecraft, Centaur, is used throughout the mission. Centaur is the upper stage for Atlas V rocket of the United Launch Alliance, with a 2,316 kg structure mass and 20,830 kg propellant capacity [39]. It uses liquid oxygen and liquid hydrogen as the propellant, with a specific impulse of 450.5 s.

**Table 2 Demand and Supply of the Resupply Mission**

| Payload Type | Node | Demand Time, day | Supply |
|---|---|---|---|
| Crew, no. | Earth | 0, 182 | 4 |
| Crew, no. | Space Station | 0, 182 | -4 |
| Science instrument, kg | Space Station | 0, 91, 182, 273 | - 1729 |
| Maintenance spares, kg | Space Station | 0, 91, 182, 273 | - 891 |
| Crew consumables, kg | Space Station | 0, 91, 182, 273 | - 1556 |
| Crew, no. | Space Station | 183, 365 | 4 |
| Crew, no. | Earth | 183, 365 | -4 |
| Science, maintenance, and crew consumables, kg | Earth | 0-365 | $+\infty$ |

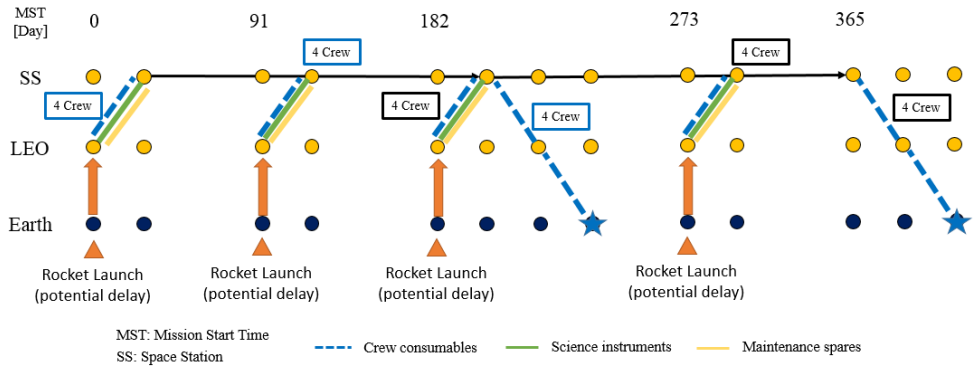

**Fig. 2 Space Station Resupply Mission Plan.**

*2. Rocket Launch Delay Scenarios*

Launch vehicle delay is one of the main uncertainty inputs for flexibility management framework. For this case study, we build a cumulative distribution function (CDF) for rocket launch delay based on the available dataset in the format of a doubly truncated exponential function, where the maximum delay is considered to be 90 days (i.e., the assumption is that if the delay is longer than 90 days, e.g., due to launch vehicle failures, we would replan the missions with a new optimization.) For a multi-mission space campaign, we assume that the rocket launch delays are independent for each mission. The argument for this assumption is that rocket launches are typically scheduled to optimize cargo transportation following specific mission requirements and there can be multiple transportation service providers for a campaign. Thus, we perform inverse transform sampling for each rocket launch mission and combine generated launch delay sampling sets together to form mission scenarios.

For result analysis, the resulting decision rules and space logistics planning are evaluated using different sets of operating mission scenarios than those used for optimization. This is to ensure the generality of the resulting decision rules. Furthermore, when we generate the Pareto front of the expected mission cost and the expected mission performance, we evaluate all solutions on the Pareto front based on a common set of evaluation mission scenarios for a fair comparison.

**B. Space Station Resupply Mission Analysis and Results**

In this section, we analyze the performance of the proposed flexibility management framework and the regular pattern of decision rules for a space station resupply campaign. For this case study, the multi-stage stochastic problem is solved under uncertainty mission scenarios generated by the method as introduced above. Each scenario contains the launch delay information of four rocket launches.

In Fig. 3, we first compare the Pareto fronts generated through different numbers of operating mission scenarios. The decision rules achieved through these mission scenario samples are evaluated based on the common evaluation sampling set with 256 mission scenarios. Results show that the sample size does not have much impact on the shape of the Pareto front. Thus, for computation efficiency, in the following analysis, we use a sample set with 128 operating mission scenarios to solve for decision rules and a different sample set with 256 evaluation scenarios to generate the Pareto front.

The decision rule strategies are also compared with traditional deterministic design strategies in Fig. 3. When the mission performance weighting coefficient $\gamma = 0$, decision-makers only focus on the mission cost of the campaign. It is equivalent to the "design for the best case" scenario, where uncertainties are neglected during mission planning and the objective is to minimize the mission cost. It is the design with the most optimistic forecast. As the value of $\gamma$ increases, mission performance becomes dominant in mission planning. The expected operating time loss decreases while the expected mission cost increases. When $\gamma$ reaches infinity, decision-makers only focus on mission performance. It is equivalent to the "design for the worse case" scenario, where mission uncertainty is the only concern during mission planning and the objective is to minimize the operating time loss. As a result, the safety stocks are prepared to counter the longest possible launch delay scenario (i.e., 90 days delay in this case study). It is the design with the most conservative forecast. The design points considering the best and the worse uncertainty scenarios are at two ends of the Pareto front. They are the anchor points for the Pareto front.

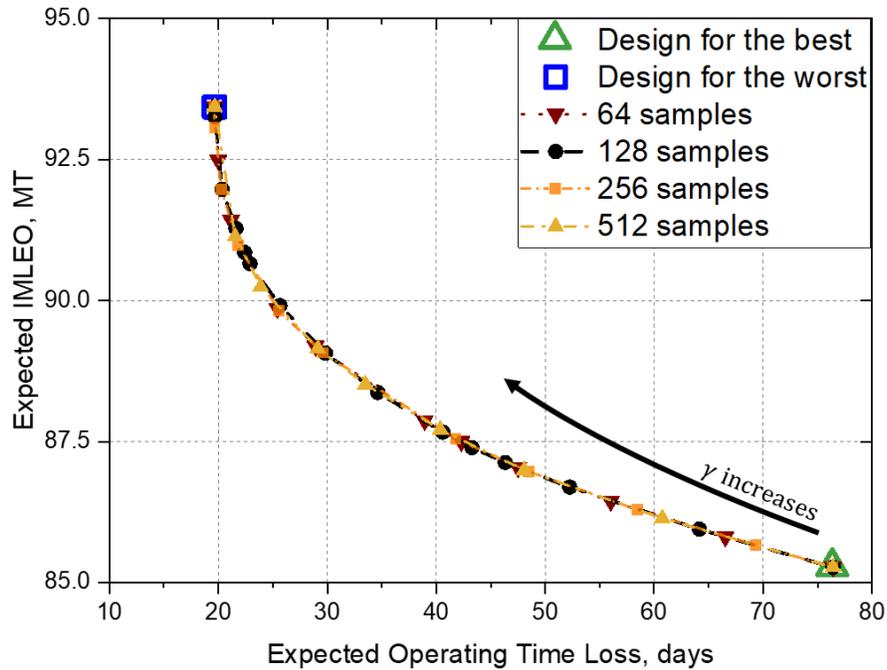

Fig. 3 Expected Operating Time Loss vs Expected IMLEO Pareto Front.

Moreover, it is also interesting to observe how the slope of the Pareto front changes as $\gamma$ increases. When the value of $\gamma$ is small, as the decrease of the operating time loss, the mission cost increases almost linearly. However, when $\gamma$ becomes large enough, the decreasing of the operating time loss leads to a much faster increase in the mission cost. To understand this slope change in the Pareto front, we conduct a sensitivity analysis to observe the variance of resultant decision rules with different values of $\gamma$, as shown in Fig. 4.

In Fig. 4, we only consider the logistics of cargo missions. One rocket launch represents one cargo transportation mission. The safety stock starts from rocket launch index 2 because we assume there is no initial safety stock available for the first rocket launch. This figure shows that when the value of $\gamma$ is small (e.g., 100), the required safety stock for each launch decreases linearly as the procedure of the space campaign. This is because the redundant supplies that are not consumed for the earlier launches can be used for later missions. Thus, more safety stock will be transported in earlier missions. However, there is an upper bound for the safety stock, where it is enough to support the longest launch delay for each mission. As the increase of $\gamma$, safety stocks for the first few launches reaches the upper bound (~2,600 kg) first. For the same amount of safety stock, it has less and less potential to counter the operating time loss. This is where the Pareto front bending begins. When $\gamma$ is above 5000, only the safety stocks for the last few launches have not been fully filled. Increasing safety stock leads to negligible improvement in mission performance. The Pareto front approaches a vertical line eventually.

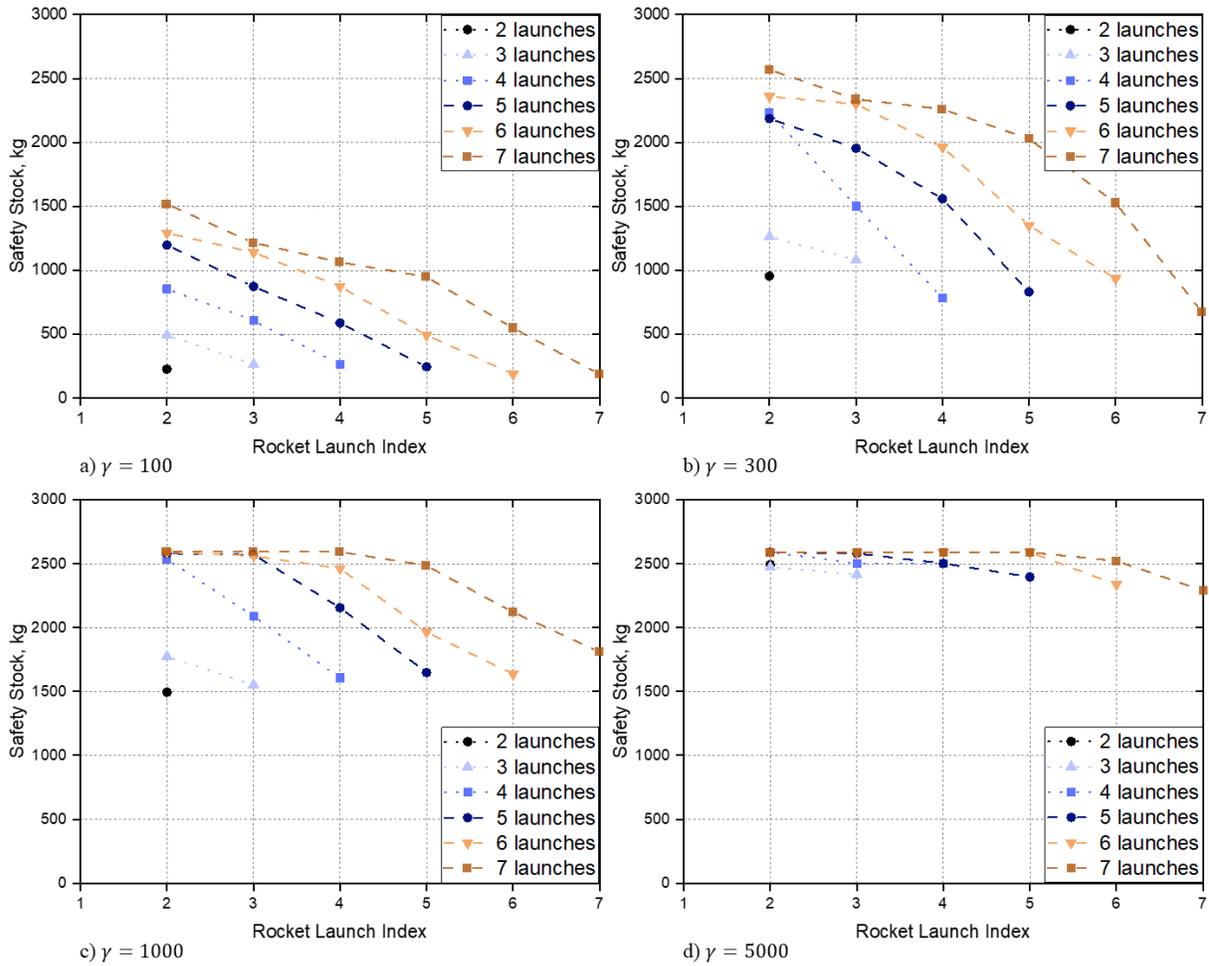

**Fig. 4 Decision Rule Result Variance.**

The impact of the number of missions is also illustrated in Fig. 4. We keep the launch frequency and mission demands but elongate the time horizon of the space campaign. In this case, we only consider cargo missions without the logistics of the crew. Longer space campaign provides larger potentials to the safety stock in the first few missions to counter the operating time loss in later missions. This means that the Pareto front bending begins later for a long-term space station operation campaign, as shown in Fig. 5. This figure shows the expected mission cost and mission performance on average per cargo mission. For different numbers of missions with similar mission demands, the slope of the Pareto front is the same when $\gamma$ is small. As $\gamma$ increases, the Pareto front for shorter space campaign with fewer space missions starts to bend first. Eventually, the Pareto fronts all end up being vertical lines. The delayed bending of Pareto front represents that decision-makers can reduce the operating time loss without a significant increase in the mission cost. Moreover, longer space campaigns also further reduce the average mission cost and improve average mission performance.

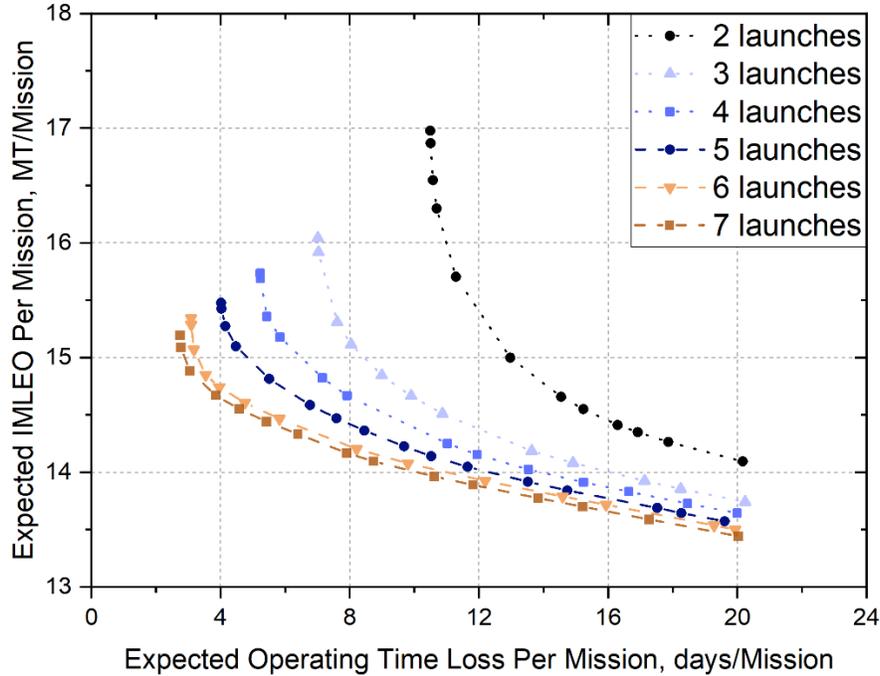

**Fig. 5 Pareto Fronts Comparison with Respect to the Number of Launches.**

Crew time assignment is another important factor that impacts the decision rules and the shape of the Pareto front. Its impact on the Pareto front is also correlated with the material consumption rate. By default, we assume 80% of the crew operating time is used for scientific experiments and 20% of the time is used for maintenance. By varying the percentage of the crew operating time used for science, we conduct a sensitivity analysis on the crew operating time assignment. The result is shown in Fig. 6. For this case study, the consumption rate for scientific experiments is twice of that for maintenance. Therefore, if 66% of the crew operating time is assigned for science, where science instruments and maintenance spares can support the same length of space station operation per unit mass, the Pareto front is the least convex, as shown in Fig. 6. On the other hand, if the material with the lower consumption rate (i.e., maintenance spares) is dominant in supporting space station daily operations, the Pareto front contains a sharp turn when balancing mission cost and mission performance. For these Pareto fronts, we can identify an area of "knee region". For multi-objective optimization, a set of solutions on the Pareto-optimal front are called "knees" when a small improvement in one objective would lead to a large deterioration in at least one other objective [40]. Focusing on the search of the knee region can generate a smaller set of solutions that are more likely to be preferred by decision-makers during mission planning. As shown in Fig. 6, outside this region, any improvement of mission performance leads to a dramatic increase in mission cost, and any decrease of mission cost leads to a dramatic deterioration of mission performance. When one of the materials is dominant in supporting space station daily operations, the knee region becomes smaller.

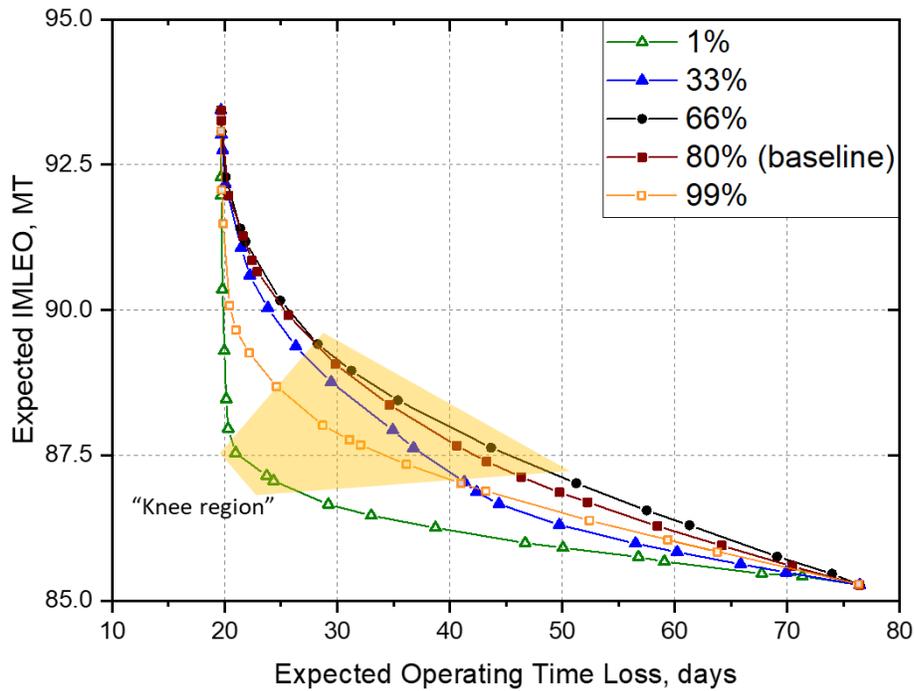

Fig. 6 Pareto Fronts for Different Percentages of Crew Operating Time for Science.

## IV. Conclusion

This paper proposes a flexibility management framework to add built-in flexibility in space logistics mission planning through decision rules and multi-stage stochastic programming. It is formulated as an optimization problem with weighted sum objectives balancing mission cost and mission performance. Space logistics optimization is developed through a network-based commodity flow model. The decision rule formulation is established by combining condition-go and linear decision rule formats to take into account the level of available safety stock before each launch and determine the amount of additional supply to be delivered together with pre-determined mission demands. A mission planning Pareto front is generated as the output to analyze the tradeoff between mission cost and mission performance.

An example mission scenario about space station resupply is established to evaluate the performance of the proposed framework considering the rocket launch delay uncertainty. Results show that the relationship among safety stocks, expected mission cost, and expected mission performance can be analyzed with the proposed framework and illustratively presented on the Pareto front between the expected mission cost (i.e., initial mass in low-Earth orbit) and the expected mission performance loss (i.e., effective crew operating time loss). Changes in space mission demands, time horizon, mission operation patterns can be directly illustrated through the variation of the Pareto front. These observations of the Pareto front show the value of the proposed flexibility management framework. For decision-makers, the decision rule results and Pareto front achieved by the framework can help improve the understanding of impacts from stochastic mission operation environments and make directly implementable strategies to counter uncertainties.

The proposed framework can be extensively applied to problems with multiple different uncertainty sources, including demand changes, spacecraft flight delays, and infrastructure performance uncertainties, and the logistics for multiple destinations at the same time. Further applications can also be explored to consider large scale human exploration and the implementation of other mitigation methods rather than only relying on the redundant supply.

## Appendix: Probability Distribution Generation for the Case Study

The launch vehicle delay distribution used in the case study is a hypothetical distribution created based on ISS past data. This Appendix includes the data source and how the used distribution is generated from the data.

The data source for the ISS U.S. On-orbit Segment (USOS) mission launch data is the ISS Flight Plan from the Flight Planning Integration Panel (FPIP) and includes ISS flights from the NASA Commercial Resupply Service (CRS) program, European Space Agency (ESA), and Japan Aerospace Exploration Agency (JAXA). The ISS Flight Plans provide a multi-year look ahead at planned mission dates for ISS activities including the CRS and international partner missions. While the complete revision history of the ISS Flight Plans is not publicly available, these plans are frequently incorporated in other NASA presentations related to ISS status and planning.

For each revision of the ISS Flight Plan at a defined baseline plan date, the planned launch for each ISS mission in the multi-year planning window is defined. These data provide insight into changes to the near and long-term missions over planning periods approaching the actual launch date. The actual launch dates for each of the ISS USOS missions were obtained from the NASA website [5]. For each ISS USOS mission, the number of days from the baseline planning date to the planned launch date was compared to the number of days the mission was delayed from the planned launch date for launches planning within 1 year of the planning, the resulting data is shown in Fig. 7. The data provided in Table 3 contains 21 launches spanning from March 2013 to February 2017 and 16 ISS Flight Plans spanning from November 2012 to January 2017 [41-43] (See Table 3). It is important to note there are more data points in Fig. 7 than ISS USOS missions because there are multiple revisions of the ISS Flight Plan between the first planning date for a given mission and the final planning date for the mission prior to launch. Additionally, the data in Fig. 7 includes the total delay between the planned mission date and the actual launch date and does not include any segregation of data for different reasons for launch delays.

**Table 3 Summary of FPIP Planning Dates**

| FPIP Dates |
|---|
| 11/13/2012 [41,42] |
| 3/26/2013 [42] |
| 7/17/2013 [43] |
| 10/15/2013 [42] |
| 11/25/2013 [42] |
| 2/6/2014 [42] |
| 4/23/2014 [42] |
| 7/3/2014 [42] |
| 8/8/2014 [42] |
| 9/25/2014 [42] |
| 12/8/2014 [42] |
| 2/20/2015 [42] |
| 10/26/2015 [42] |
| 2/18/2016 [42] |
| 9/20/2016 [42] |
| 1/20/2017 [42] |

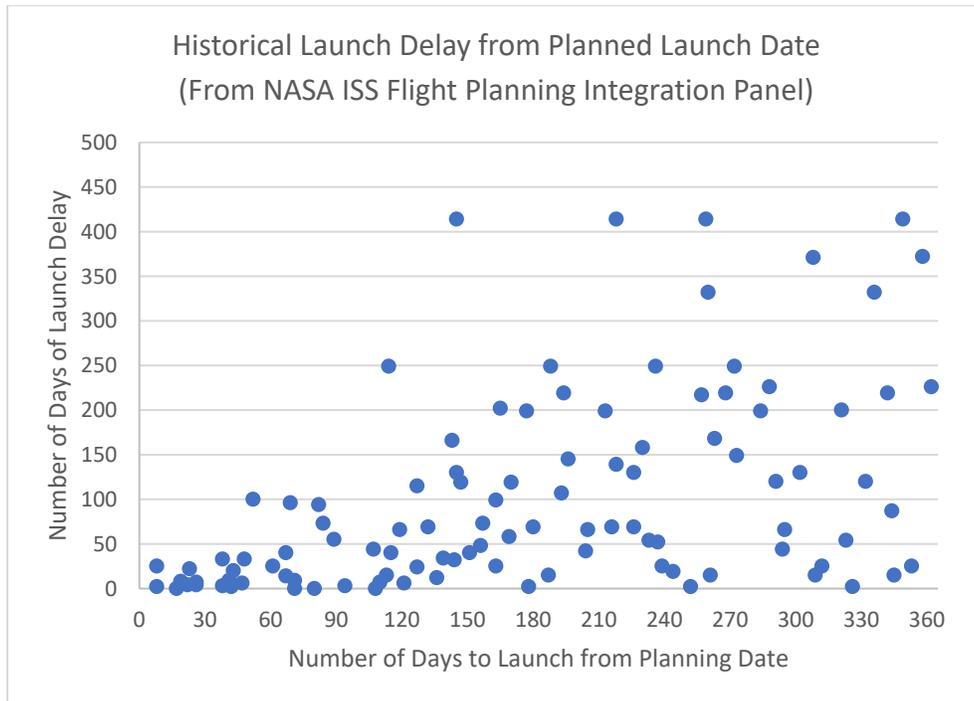

**Fig. 7 Historical Error in CRS Launch Date Planning (within 1 Year of Planned Launch).**

The data in Fig. 7 shows an expected trend that as the mission planning date is closer to the planned launch date, the accuracy of the planned launch date improves. Additionally, variation in the launch delay increases as the number of days between the mission planning date and the launch date increases. The data points with delays exceeding 350 days are associated with return to flight missions for Orbital ATK and SpaceX after the Orb-3 and SpX-7 launch failures, respectively, and the HTV4 and HTV5 missions from JAXA.

For a multi-mission space station resupply campaign, based on the available data, the launch delay for each mission is assumed to be independent in the scenario generation. There are two reasons for this assumption. Firstly, the ISS Flight Plan provides mission plans based on planned ISS needs well in advance of the actual mission date. As the ISS planning is adjusted based on ISS needs, the mission dates may shift to optimize cargo delivery to maximize total delivered mass, providing the best value for the NASA and the taxpayer. Secondly, the CRS program has two launch providers, so if one provider experiences a significant launch delay, NASA can reorganize launches to minimize delay in providing cargo to the ISS.

Because the launch delays are assumed as independent for each mission, we only need to discuss the launch delay model for one mission first and implement it directly to other missions in the campaign. In this case study, we assume four missions per year as a campaign; so we only use the data within a 90-day planning period for the launch vehicle delay model. For simplicity in this case study, we consider 90 days as the longest possible delay; we assume that if the delay is longer than 90 days, we would replan the missions with a new optimization. The data for the 90-day planning window is shown in Fig. 8. The cumulative density function of launch delay (in the number of days) within a 90-day planning window is determined from the data in Fig. 8. The resulting discrete probability is curve fit to find a continuous function to represent the launch delay. The resulting curve fit is a doubly, truncated exponential curve shown in Fig. 9. Based on the launch delay cumulative density function as shown in Fig. 9, we can do inverse transform sampling to generate uncertainty scenarios input for the space station resupply mission example.

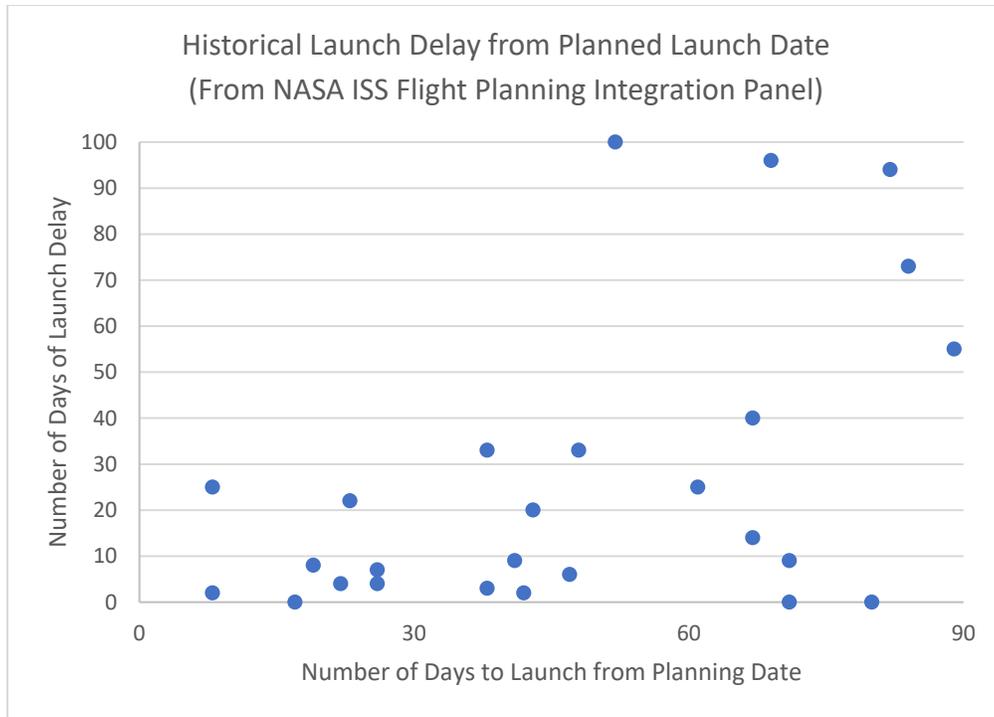

**Fig. 8 Historical Error in CRS Launch Date Planning (within 90 days of planned launch).**

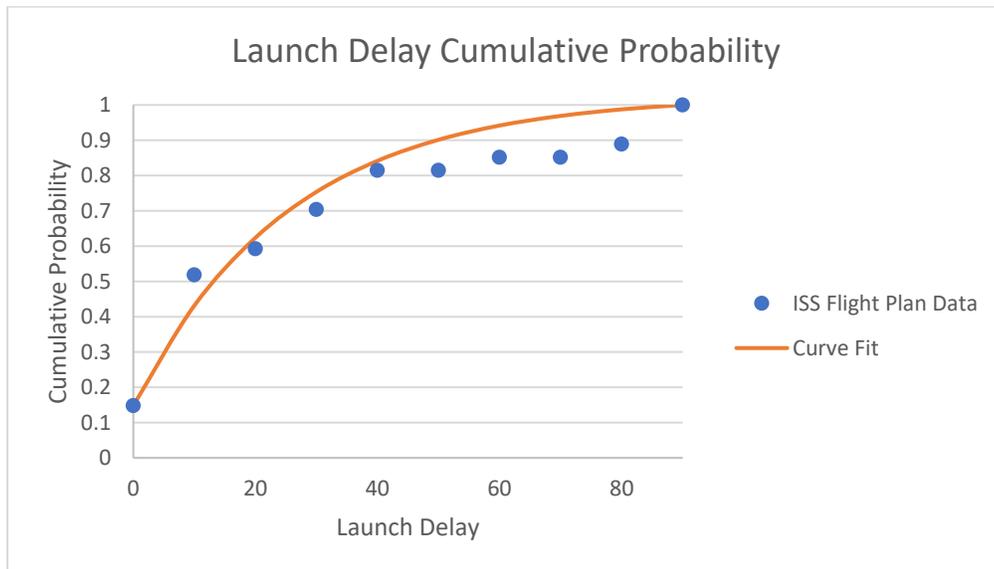

**Fig. 9 Launch Delay Cumulative Distribution Function.**